\documentclass[12pt]{article}
\usepackage{amssymb}
\usepackage{amsmath, amscd, amsthm}
\usepackage{mathrsfs}
\usepackage{amsfonts}
\usepackage{amsrefs}
\begin{document}

\theoremstyle{remark}
\newtheorem*{Remark}{Remark}
\newtheorem*{Remarks}{Remarks}
\newtheorem*{Proof}{Proof}
\theoremstyle{definition}
\newtheorem*{Definition}{Definition}
\newtheorem*{Example}{Example}
\newtheorem*{Examples}{Examples}

\theoremstyle{plain}
\newtheorem*{Theorem}{Theorem}
\newtheorem*{Theorema}{Theorem}
\newtheorem*{Proposition}{Proposition}
\newtheorem*{Corollary}{Corollary}
\newtheorem*{Lemma}{Lemma}

\newcommand{\g}{\mathfrak g}
\newcommand{\cl}{\mathbf{Cl}}
\newcommand{\vs}{\vspace{5mm}}
\renewcommand{\labelenumi}{\roman{enumi})}
\setlength{\parindent}{0em}
\small\normalsize
\begin{center}
{\Large\sc Analytic Tate spaces and reciprocity laws\footnote{\today}}\\ $\ $ \\
{by Ricardo Garc\'{\i}a L\'opez\footnote{
Departament d'\`Algebra i Geometria.
Universitat de Barcelona,
Gran Via, 585.
E-08007, Barcelona, Spain.
e-mail: ricardogarcia@ub.edu}
}
\end{center}\vs
\begin{abstract}
We consider an analytic variant of the notion of Tate (or locally linearly compact) space 
and we show that, both in the complex and in the $p$-adic analytic setting, one can use it to 
define symbols which satisfy Weil-type reciprocity laws for curves.
\end{abstract}
\vs

{\bf 1. Introduction}\vs 

A topological vector space is said to be linearly compact if it is the topological dual of a vector space endowed with the discrete topology.
A Tate space is a topological vector space which contains a linearly compact open subspace 
(according to \cite{Dri}, this terminology is due to A. Beilinson. They were considered by S. Lefschetz 
under the name of locally linearly compact spaces, see \cite{Kap} and the references therein). Tate spaces play 
a relevant r\^ole  in the study of algebraic curves and chiral algebras (see, among others, [loci cit.],  \cite{BBE}, \cite{BD}).\vs

A Tate space splits as the topological direct sum of a discrete and a linearly compact space. The simplest example is 
the ring $k((t))$ of Laurent power series with coefficients in a field endowed with the $t$-adic topology, 
one has the topological splitting 
$k((t))=t^{-1}k[t^{-1}]\oplus k[[t]]$ where the first summand is discrete and the second is linearly compact. 
The analogous decomposition for convergent power series suggests that an analytic counterpart of Tate spaces might be
the category of those topological vector spaces which split as the sum of a nuclear Fr\'echet space  
and the strong dual of a nuclear Fr\'echet. Examples of such spaces are
the space of germs of analytic functions on a punctured neighborhood of a point in the complex plane 
or the Robba rings appearing in the theory of $p$-adic differential equations. \vs

Our aim in this note is to show that some results related to the notion of Tate space are still valid in the analytic
setting if one considers the category described above. In section 2 we collect some results 
from functional analysis, in section 3 we define commutator symbols in this context, and
in section 4 we prove reciprocity laws for curves, both in the complex and $p$-adic analytic cases, following similar arguments as those used in the formal case. \vs

I thank the anonimous referee for her/his comments. \vs

{\bf 2. Analytic Tate spaces}\vs

Let $k$ denote a local field of characteristic zero, set $k^{\times}=k-\{0\}$. In what follows, we refer 
to \cite{Tre}, \cite{MV} (in the archimedean case) and to \cite{Sch} (in the non-archimedean case) for unexplained terminology.
All topological $k$-vector spaces will be assumed to
be locally convex, direct sums will be assumed to have the locally convex direct 
sum topology. The  dual of a vector space $V$ will be denoted $V^{\ast}$ 
(in case $V$ is a topological vector space, this notation will refer to the strong dual), the dual of a 
map $g$ will be denoted $g^{\ast}$. 
A (FN)-space is a vector space which is simultaneously Fr\'echet and nuclear. A (DFN)-space is the 
strong dual of a (FN)-space. The following theorem summarizes some results to be used later on:\vs

(2.1) {\bf Theorem}:{\it \,
\begin{enumerate}
\item A topological vector space $V$ is a (FN)-space (respectively, a (DFN)-space) 
if and only if there is a sequence of locally convex topological
vector spaces and nuclear maps (resp., nuclear injective maps)
\begin{eqnarray*}
&V_1\longleftarrow V_2& \longleftarrow\dots \\
(&V_1\longrightarrow V_2& \longrightarrow\dots)
\end{eqnarray*}
such that $V=\varprojlim V_i$ (resp., $V=\varinjlim V_i$).
\item If a space is simultaneously (FN) and (DFN), then it is finite dimensional.
\item Closed subspaces and separated quotients of a (FN)-space
are (FN)-spaces. 
\item Closed subspaces and separated quotients of a
(DFN)-space are (DFN)-spaces.
\end{enumerate}}

Here, i) is a small variation on the results 
of \cite{Kom} and can be proved in the same way, and ii), iii), iv) are well-known 
(see e.g \cite{FW}*{\S 25, 2.11}, \cite{Cr}*{Proposition 2.8}), 
with the exception of the assertion about closed subspaces of (DFN)-spaces in the non-archimedean case, 
which is proved in \cite{Cr}*{Lemma 1.2}. \vspace{1mm}

One can transpose to the analytic setting the notion of Tate space
(see \cite{Dri}*{3.1}) as follows:\vspace{3mm}

\begin{Definition}
A topological vector space $V$ will be said to
be an \emph{analytic Tate space} if there exist a (FN)-space $L\subset V$
and a (DFN)-space $G\subset V$ such that  $V= L \oplus G$ topologically .
In this case, we
will say that $L$ is a (FN)-lattice and $G$ is a (DFN)-lattice in
$V$. By a lattice we will refer indistinctly to a (FN) or a (DFN)-lattice.
\end{Definition}\vspace{1mm}

\begin{Remark}
To some extent, lattices  are
the analogues of the compact and discrete lattices considered in \cite{Dri}.
For example, as in the locally linearly compact case,
if $L_2\subset L_1$ are lattices, then $L_1/L_2$ is finitely dimensional  (being a closed subspace of $V/L_2$ and a separated quotient 
of $L_1$, it is simultaneously
a (FN) and a (DFN)-space). 
However, contrarily to what happens in the formal case, in general neither the intersection nor the sum of two 
lattices is a lattice, see the remark after example i) below. \vs

(2.2) {\bf Examples}: i) (see \cite{ADKP}*{$\S$ 1}) Let $k=\mathbb C$. For $U\subset \mathbb C$ open, 
denote $H(U)$ the ring of analytic 
functions on $U$, let 
\[
\mathscr O_{\widehat 0}=\varinjlim_{0\in U} H(U-\{0\})
\]
denote the ring of germs 
of analytic functions on a punctured neighborhood of $0\in\mathbb C$, put 
\begin{eqnarray*}
 \mathscr O_{\widehat 0}^+&=&\varinjlim_{0\in U} H(U)\\
 \mathscr O_{\widehat 0}^-&=&\varprojlim_{0\in U} 
H_{\infty\mapsto 0}(\mathbb C - \overline{U})\ ,
\end{eqnarray*}
where $U$ runs over a neighborhood basis of $0$ and the notation $\infty\mapsto 0$ means that we consider only those  functions $f$
such that $\lim_{\|x\|\mapsto \infty}f(x)=0$. Then $ \mathscr O_{\widehat 0}^-$ is a (FN)-space,
$\mathscr O_{\widehat 0}^+$ is a (DFN)-space and we have a topological isomorphism
\[
 \mathscr O_{\widehat 0} \cong \mathscr O_{\widehat 0}^+  \oplus \mathscr O_{\widehat 0}^- \ ,
\]
thus $ \mathscr O_{\widehat 0} $ is an analytic Tate space. Once a coordinate has been fixed, 
one can construct similar examples considering 
power series with Gevrey conditions (in particular, the $1$-Gevrey case can be seen as 
the ring of microdifferential operators with constant coefficients). \vs

\begin{Remark}
If $f$ is a unit of the ring $\mathscr O_{\widehat{0}}$, then
we have two decompositions
\[
\mathscr O_{\widehat{0}}=\mathscr O^+_{\widehat{0}} \oplus \mathscr O^-_{\widehat{0}}=
f\cdot \mathscr O^+_{\widehat{0}} \oplus f\cdot \mathscr O^-_{\widehat{0}}\, ,
\]
so $\mathscr O^+_{\widehat{0}}$ and $f\cdot \mathscr O^+_{\widehat{0}}$ are two 
(DFN)-lattices in $\mathscr O_{\widehat{0}}$. If $f$ has an essential singularity at zero, then
\[
\mathscr O^+_{\widehat{0}} \cap (f\cdot \mathscr O^+_{\widehat{0}})=\{ 0\} ,
\]
thus an intersection of lattices does not need to be a lattice. Also, if 
$\mathscr O^+_{\widehat{0}} + f\cdot \mathscr O^+_{\widehat{0}}
$ were a (DFN)-lattice, there would exist a (FN)-subspace $F\subset \mathscr O_{\widehat{0}}$
with $\mathscr O^+_{\widehat{0}} \oplus (f\cdot \mathscr O^+_{\widehat{0}}) \oplus F = \mathscr O_{\widehat{0}}$.
But then
\[
 (f\cdot \mathscr O^+_{\widehat{0}}) \oplus F \cong 
\mathscr O_{\widehat{0}}/\mathscr O^+_{\widehat{0}} \cong \mathscr O^-_{\widehat{0}}
\]
which is a (FN)-space, so that $\mathscr O^+_{\widehat{0}} \cong f\cdot \mathscr O^+_{\widehat{0}}$ would be simultaneously 
(FN) and (DFN), which is impossible (by 2.1.ii) . It follows that, in general, sums of lattices are
not lattices either.
\end{Remark}
\vs

ii) The ring $k[[t]]$, endowed with the topology given 
the bijection $k[[t]]\cong \prod_{\mathbb N}k$
is a (FN)-space (see e.g.\cite{Tre}). Its strong dual can be identified with a polynomial ring $k[u]$, 
endowed with the topology given by the bijection $k[u]\cong \oplus_{\mathbb N}k$.
Thus, the topological direct sum $k((t))=k[[t]] \oplus t^{-1}k[t^{-1}]$ is an analytic Tate space.
\vs

iii) If $\mathfrak g$ is a finite dimensional topological $k$-Lie algebra, then the loop algebra
$\mathfrak g ((t)) := \mathfrak g \otimes_k k((t))$ is an analytic Tate vector space and the
Lie braket
\[
 [g\otimes p , h\otimes q]= [g,h]\otimes p\,q
\]
is continuous, thus $\mathfrak g ((t))$ is an analytic Tate Lie algebra. It is likely that one can give a construction of
semi-infinite cohomology in this setting, following the method in \cite{BD}*{3.8}. \vs

iv) (see \cite{Cr0}*{Part II}, \cite{Cr}) Assume $k$ is non-archimedean. 
Given an interval $I\subset [0, +\infty)$,
let $A(I)$ be the ring of power series $\sum_{i\in\mathbb Z}\,
a_i\,t^i$, $a_i\in k$, which are convergent for $\lvert t\rvert \in
I$. For each $\rho\in I$, the ring $A(I)$ is endowed with the
norm
\[
\left\|\ \sum_{i\in\mathbb Z}\, a_i\,t^i \ \right\|_{\rho} = \sup_{i\in\mathbb
Z}\mid a_i\mid \rho^{i}
\]
and, with the topology defined by this family of norms, the
$k$-algebra $A(I)$ is a Fr\'echet space. The union
\begin{equation*}
\mathcal R= \bigcup_{r<1}A\,(r,1)=\left\{\sum_{i\in\mathbb Z} a_i\,t^i \ \middle\vert \ a_i\in k,  \left.
\begin{array}{rl}
 \exists \lambda <1 & \lvert a_i\rvert\lambda^i \rightarrow 0 \mbox{ for } i \to -\infty \\
\forall \varepsilon <1 & \lvert a_i\rvert\varepsilon^i \rightarrow 0 \mbox{ for } i \to +\infty \\
\end{array}\right.
\right\}
\end{equation*}
is a ring, called the Robba ring (over the field $k$). One considers in $\mathcal R$ the
direct limit topology given by the first equality above. If we put $\mathcal
R^{+}=\mathcal R\cap
 k[[t]]$,
$\mathcal R^{-}=\mathcal R\cap k[[t^{-1}]]$, and we endow these
spaces with the subspace topology, then the strong dual of
$\mathcal R^{+}$ is isomorphic to $t^{-1} \mathcal R^{-}$ via the residue pairing
\[
\left\langle\ \sum_{i\geqslant 0}a_i\,t^i\,,\sum_{i\geqslant
0}b_i\,t^{-i-1}\ \right\rangle=\sum_{i\geqslant 0}a_i\,b_i\,,
\]
and the direct sum decomposition $\mathcal R=\mathcal
R^{+}\oplus t^{-1}\mathcal R^{-}$  is topological. 
It is not difficult to prove that
$\mathcal R^{+}$ is a (FN)-space, so $t^{-1} \mathcal R^{-}$ is a (DFN)-space and
the Robba ring $\mathcal R$ is an analytic Tate space (notice  the properties of the spaces of positive and 
negative powers of $t$ are reversed with respect to those in example i), namely there $\mathscr O^+_{\widehat{0}}$
is (DFN) and $\mathscr O^-_{\widehat{0}}$ is (FN)). See \cite{Cr0}*{7.2} for a coordinate free description of the Robba ring attached 
to a point of a smooth curve defined over a finite field. \vs
\vs

{\bf 3. Polarizations and pairings} \vs

 We recall the main features of the theory of  polarizations and determinants of Fredholm maps
(see \cite{Seg}*{Appendix B}, cf. \cite{AR}, \cite{Bry}, \cite{Kap}):\vs

Let $Pic^{\mathbb Z}$ denote the category of graded lines. An object in
$Pic^{\mathbb Z}$ is a pair $(\ell, n)$ where $\ell$ is a $1$-dimensional $k$-vector space and
$n\in\mathbb Z$. In this category, $Hom((\ell_1,n_1), (\ell_2,n_2))=\emptyset$ if $n_1\neq n_2$
and $Hom((\ell_1,n), (\ell_2,n))=Hom_k(\ell_1, \ell_1)$. 
There is a tensor product $\otimes: Pic^{\mathbb Z}\times  Pic^{\mathbb Z} \longrightarrow Pic^{\mathbb Z}$ defined by
\[
((\ell_1,n_1), (\ell_2,n_2)) \longmapsto (\ell_1\otimes \ell_2, n_1+n_2)
\]
And a commutativity constraint
\[
(\ell_1\otimes \ell_2, n_1+n_2) \stackrel{\sim}\longrightarrow (\ell_2\otimes \ell_1, n_2+n_1)
\]
given by $u\otimes w \longmapsto (-1)^{n_1n_2} w\otimes u$. If $(\ell, n)$ is a graded line, 
the reference to the integer $n$ will be omitted when no confusion arises. \vs

Let $V, W$ be topological vector spaces. A continuous linear map $f: V \longrightarrow W$ is said to be Fredholm if it has finite dimensional kernel
and cokernel. In this case, set $n=\dim \mbox{Ker}(f)$, $m=\dim\mbox{Coker}(f)$, 
$\mbox{index}(f)=m-n$, and put
\[
 Det(V,W,f)= \left( \left(\bigwedge^{n}\mbox{Ker}(f)\right)^{\ast} \otimes \bigwedge^{m}\mbox{Coker}(f) \, , \
\mbox{index}(f) \,\right).
\]

Assume $V,W$ are (FN)-spaces and $f:  V \longrightarrow W$ is Fredholm. 
If $u: V \longrightarrow W$ is nuclear then,  $\mbox{index}(f)=\mbox{index}(f+u)$ (\cite{CM}*{8.1.-1.}).
This statement holds also for (DFN) spaces, because on the one hand the strong dual of a (DFN) space is (FN), in particular  reflexive (see \cite{Sch}*{III.16.10, IV.19.3}, \cite{Kom}*{Theorem 12}), and on the other hand if $f$ is a Fredholm map between Fr\'echet spaces then the dual map $f^{\ast}$ is also Fredholm and $\mbox{index}(f)=-\mbox{index}(f^{\ast})$ (for a local field, the proof in the non-archimedean case is the same as in the archimedean case, see \cite{CM-0}). In both cases there is a canonical isomorphism
\[
 Det(V,W,f) \longrightarrow   Det(V,W,f+u)\,,
\]
this follows from the comparison between the above definition of the determinant line and the one 
given by G. Segal in \cite{Seg}, which is invariant under nuclear perturbations (see \cite{Seg}*{Appendix B}, 
see also \cite{Hu}*{D.2.11}). We put 
\[
\mathbf{Det}(V,W,f)=\varinjlim_u Det(V,W,f+u)
\]
where the transition morphisms are the isomorphisms above.  \vs

(3.1) {\bf Definition}: On the set of topological direct sum 
decompositions $V=V_+\oplus V_-$ we consider the following equivalence relation:
$
V_+\oplus V_- \sim W_+\oplus W_- 
$
if and only if for $i,j\in\{+,-\}$, $i\neq j$ the compositions 
\[
V_i \hookrightarrow V \twoheadrightarrow W_j \ \ \mbox{ and } \ \
W_i \hookrightarrow V \twoheadrightarrow V_j
\]
are nuclear, where in both cases the arrows $\twoheadrightarrow$ 
are the projections attached to the given decompositions. A polarization of $V$ is an equivalence
class of decompositions. If we fix a decomposition $V=V_+\oplus V_-$, the decompositions in 
the same equivalence class will be called allowable, as well as the projections onto its summands. 
\vs

Let $V_{+}, W_{+}$ be two allowable plus-summands, 
let $p: V \longrightarrow W_+$ be an allowable projection. The restriction of $p$ to $V_+$ 
is Fredholm (since it is invertible modulo a nuclear operator). 
If $p': V \longrightarrow W_{+}$ is another allowable projection, the difference 
$p_{\mid V_+} - p'_{\mid V_+}$ is nuclear, so we have a canonical isomorphism
\[
  \mathbf{Det}(V_{+},W_{+}, p_{\mid V_+}) \cong \mathbf{Det}(V_{+},W_{+},p'_{\mid V_+})\
\]
and we put 
\[  
\mathbf{Det}(V_{+} : W_{+})=\varinjlim_{p} \mathbf{Det}(V_{+},W_{+},p_{\mid V_+})\,,
\]
where the limit runs over all restrictions to $V_+$ of allowable projections $V \twoheadrightarrow W_{+}$. 
If $U_+,V_+,W_+$ are allowable, the composition of two allowable projections
$V_+ \twoheadrightarrow U_+ \hookrightarrow V \twoheadrightarrow W_+$ differs from an allowable projection
$V_+ \twoheadrightarrow W_+$ by a nuclear perturbation, and so it follows that there are canonical isomorphisms 

\[
\omega(V_+,U_+,W_+): \mathbf{Det}(V_{+} : W_{+}) \cong \mathbf{Det}(V_{+} : U_{+}) \otimes \mathbf{Det}(U_{+} : W_{+})\, .
\]

Let $V$ be a polarized vector space, choose an allowable decomposition $V=V_+\oplus V_{-}$,
denote $GL(V)$ the group of bicontinuous automorphisms of $V$.
Given $g\in GL(V)$, denote $g_{i,j}:V_i\longrightarrow V_j$ the composition of the restriction $g_{\mid V_i}$ 
with the allowable projection $V \twoheadrightarrow V_j$, where $i,j\in\{+,-\}$. The restricted linear group of $V$ is  the group
\[
 GL_{res}(V) = \left\{ g  \in GL(V) \ \mid \  
g_{+,-} \ \mbox{ and }\  g_{-,+} \mbox{ are nuclear } \right\}
\]
Notice that if $g_{+,-}, g_{-,+}$ are nuclear for one allowable decomposition, then they are 
so for all of them and $(g^{-1})_{+,-}, (g^{-1})_{-,+}$ are also nuclear. \vs

Choose allowable subspaces $V_+, W_+\subset V$. For $g\in GL_{res}(V)$, the subspaces $g(V_+), g(W_+)$ are 
also allowable and there is a canonical isomorphism
\[
\mathbf{Det}(V_+ : g(V_+))\longrightarrow \mathbf{Det}(W_+ : g(W_+))\, ,
\]
we denote
\[
\mathbf{P}_g=\varinjlim_{V_+} \, \mathbf{Det}(V_{+}:g(V_{+}))\, .
\]

Given $f,g\in  GL_{res}(V)$ and an allowable $V_+$, the subspaces $g(V_+), fg(V_+)$ are also allowable, 
let $\tau_f: \mathbf{Det}(V_{+}:g(V_{+})) \longrightarrow \mathbf{Det}(f(V_{+}):f(g(V_{+})))$ be the conjugation isomorphism 
which sends $\alpha: V_+ \longrightarrow g(V_+)$ to $f \circ \alpha\circ f^{-1}: f(V_{+}) \longrightarrow f(g(V_{+}))$. We have a map
\[
\mathbf{Det}(V_{+}:f(V_{+}))  \otimes \mathbf{Det}(V_{+}:g(V_{+})) \longrightarrow \mathbf{Det}(V_{+}:fg(V_{+})) 
\]
defined as the composition $\omega(V_+,f(V_+), fg(V_+))\circ (id\otimes \tau_f)$, which induces an isomorphism
\[
 \rho_{f,g}:\mathbf{P}_f\otimes \mathbf{P}_g \longrightarrow \mathbf{P}_{fg}\ .
\]
(3.4) Put
\[
 GL^+_{res}(V)^{}=\left\{(f,\alpha) \mid \ f\in  GL_{res}(V) \ ,\, \alpha\in \mathbf{P}_f \ \right\}\, ,
\]
and consider the operation
\[
(f,\alpha) \cdot (g, \beta) = (f\cdot g, \ \rho_{f,g}(\alpha\otimes \beta)) \, .
\]
With this operation, $ GL^+_{res}(V)^{}$ is a group and we have a central extension
\[
 1 \longrightarrow k^{\times} \longrightarrow  GL^+_{res}(V)^{} \longrightarrow  GL_{res}(V) \longrightarrow 1
\]
where the first map is $\lambda \mapsto (Id, \lambda)$ (notice that $\mathbf{P}_{Id}=k$ canonically), and the second is the projection. 
It is well-known that such an extension defines a symbol: Given commuting elements $f,g\in GL_{res}(V)$, 
choose liftings $\tilde{f}, \tilde{g} \in GL^+_{res}(V)$
and define the symbol $(f,g)_{+}$ by
\[
(f,g)_{+}=\tilde{f}\,\tilde{g}\,\tilde{f}^{-1}\tilde{g}^{-1} \in k^{\times}.
\]
Then $(f,g)_{+}$ is independent of the chosen liftings and it follows from the definitions that one has:
\begin{enumerate}
\item $(f,g)_{+}=(g,f)_+^{-1}$ for every commuting $f,g\in GL_{res}(V)$.
\item $(f_1\cdot f_2,g)_+=(f_1,g)_+\cdot(f_2,g)_+$ and $(f,g_1\cdot g_2)_+=(f,g_1)_+\cdot(f,g_2)_+$ 
under the corresponding commutativity assumptions.
\end{enumerate}
\

{\it Remarks}:
i) Notice that the definition of the commutator pairing given in \cite{ADK} or \cite{AR} cannot be directly transposed here
(see the remark after example i) above).\vspace{3mm}

ii) We can interchange the r\^oles of plus and minus summands,
in this way for $g\in GL_{res}(V)$ we 
obtain minus lines $\mathbf M_g=\varinjlim_{V_-} \, \mathbf{Det}(V_{-}:g(V_{-}))\, $ and a central extension
\[
 1 \longrightarrow k^{\times} \longrightarrow  GL^-_{res}(V)^{} \longrightarrow  GL_{res}(V) \longrightarrow 1\ .
\]
Again as above, we have a pairing $(\cdot, \cdot)_-$
defined for any pair of commuting elements of $GL_{res}(V)$. In the formal case, 
this pairing does not provide any new information because $(f,g)_+=(f,g)_-^{-1}$ (see \cite{AR}*{Proposition 3.3.4}).
We do not know if this is true in the present case as well, although explicit computation in the examples 
considered in the previous section shows that it does hold in these cases. \vs

The following proposition is probably well-known to functional analysts, although we have not found a precise reference. 
In any case, it is the crucial point 
for proving the results of this note.\vs

(3.5) {\bf Proposition}: {\it Let $V$ be a (FN)-space and $W$ a (DFN)-space. Then:
\begin{enumerate}
\item Every continuous linear map $f:V \longrightarrow W$ is nuclear.
\item Every continuous linear map $g: W \longrightarrow V$ is nuclear.
\end{enumerate}}

\begin{Proof} i) In the archimedean case, we can assume that there is an inductive system of Banach spaces $\{W_i\}_{i\in\mathbb N}$ with injective transition maps such that $W=\bigcup_{i\in \mathbb N} W_i$ (see \cite{MV}*{Proposition 25.20}). By Grothendieck's factorization theorem (\cite{MV}*{Theorem 24.33}) there is 
an $i\in\mathbb N$ such that $f$ factors as 
$V\longrightarrow W_i \hookrightarrow W$. Since $V$ is nuclear, $V\longrightarrow W_i$ is a nuclear map,
and then so is $f$. The same proof works in the non-archimedean case, using \cite{Sch}*{I.8.9}.  \vs

ii) Consider first the archimedean case. Since $W$ is a (DFN)-space, it is complete 
(\cite{Tre}*{Cor. 2., pg. 344}) and nuclear (\cite{Tre}*{Proposition 50.6}).
It is also barreled because it is  a 
direct limit of barreled spaces. 
Then, statement ii) is a consequence of \cite{DC}*{Theorem (1.3)}. \vs

In the non-archimedean case, it also holds that (DFN)-spaces are complete (\cite{Sch}*{III.16.10.i}),
nuclear (\cite{Sch}*{IV, 19.9}) and barreled. Inspection of the proof of \cite{DC}*{Theorem (1.3)} shows that it is
also valid in this case (at least for a spherically complete base field, which covers our setting), 
the results needed for the proof are \cite{Sch}*{IV.18.8}, and the Banach-Steinhaus theorem \cite{Sch}*{I.6.15}.
$\Box$
\end{Proof}
\vspace{2mm}

An analytic Tate space structure on a vector space $V$ defines a polarization in an obvious way. 
For polarizations defined in this manner,
it follows from the previous proposition that we have:\vs

(3.6) \, {\bf Corollary}: {\it
\, i) An analytic Tate space has an essentially unique polarization, that is, if we have a topological isomorphism
\[
V_+ \oplus V_- \cong W_+ \oplus W_- 
\]
where $V_+, W_+$ are (FN)-spaces and $V_-,W_-$ are (DFN)-spaces, then these decompositions are equivalent 
in the sense of definition (3.1).\vs

ii) If $V$ is an analytic Tate space, then $GL(V)=GL_{res}(V)$. }

\vs

{\it Remarks}: 
i) For the analytic Tate space in example (2.2.i) above, proposition (3.5) is implicitely proved in \cite{ADKP}.\vspace{3mm}

ii) The situation is somehow more symmetric than in the formal case. There,
a linear continuous map $C \longrightarrow D$ from a linearly compact to a discrete space has finite dimensional image, but
this does not hold for morphisms $D \longrightarrow C$ (e.g., for $k[t] \hookrightarrow k[[t]]$).\vspace{3mm}
 
iii) Kapranov's construction of measures and Fourier transforms in the locally 
linearly compact case (see \cite{Kap}) can be mimicked in the present context. The main point is that  if $L_1\subset L_2$ 
are lattices,  the quotient 
$L_1/L_2$ is finite dimensional, and this is the essential fact which is needed to reproduce Kapranov's constructions. For example, if $V$ is an analytic Tate space and $V_+\subset V$ is a fixed 
(FN)-lattice (respectively, (DFN)-lattice), then the assignment
\[
 U \longmapsto \mbox{Det}(V_+,U)
\]
defined on the grassmannian of (FN)-lattices (resp., (DFN)-lattices) defines a determinant theory in the 
sense of loc. cit.. Similarly, it is easy to adapt to analytic Tate spaces the definition given in \cite{Kap}*{3.2.2} 
of the semi-infinite de Rham complex associated to a locally linearly compact $\mathbb R$-vector space. 
\end{Remark}

From the previous proposition we can also derive the following finiteness result:

\begin{Corollary}
 Let $V$ be an analytic Tate space. If $L$ is a (FN)-lattice and $G$ is a (DFN)-lattice,
then $L\cap G$ is of finite dimension and $L+G$ is of finite codimension in $V$.
\end{Corollary}

 The first statement follows directly from (2.1.ii).
For the second, 
we first recall the following lemma of H. Lang (\cite{La}*{Lemma
2.2}): Let $C_1, C_2$ be two closed subspaces of a locally
convex vector space $V$. Assume $C_1$ is a topological direct
summand of $V$, let $\pi: V \longrightarrow C_1$ be a continuous
projection. If $\pi_{\mid C_2}$ is compact, then $C_1+C_2$ is
closed in $V$. \vs

It follows from this lemma and (3.5) that $L+G$ is closed in $V$. Now $V/L+G$ is the quotient 
of the (DFN)-space $V/L$ by the closed subspace
$L+G/L$, and thus is a (DFN)-space. It is also the quotient of the (FN)-space
$V/G$ by $L+G/G$. Thus, by (2.1.ii), it
is finite dimensional. \vs
$\Box$

{\bf 4. Calculations and reciprocity laws} \vs

Let $\mathscr O_{\widehat 0}$ be as in example i) above, let
$\mathscr O_{\widehat 0}^{\ast} \subset \mathscr O_{\widehat 0}$ be the multiplicative subgroup of invertible germs. 
Given $f\in \mathscr O_{\widehat 0}^{\ast}$, we 
denote also by $f\in GL(\mathscr O_{\widehat 0})$ the operator of multiplication by 
$f$ in $\mathscr O_{\widehat 0}$.
If we fix a coordinate $t$, any $f\in \mathscr O_{\widehat 0}^{\ast}$ can be written as a product
\[
 f= c\cdot t^n \cdot g(t) \cdot h(t^{-1})\ ,
\]
where $c\in \mathbb C$, $n\in\mathbb Z$, $g$ is holomorphic at zero, $g(0)=1$, $h(x)=e^{\varphi(x)}$ and $\varphi$
is an entire function with $\varphi(0)=0$
(Weierstrass-Birkhoff decomposition, see \cite{Bir}).  \vs

In the $p$-adic case (example (2.2.iv)),
because of results of E. Motzkin and G. Christol (see \cite{Chr}),
the Weierstrass-Birkhoff decomposition holds in the following form: 
Let $f\in \mathcal R^{\ast}$ be a unit of the Robba ring.
Then, there exist unique $n\in\mathbb Z$, $c\in k$, 
$g(t)\in1 +t\mathcal R^+$ invertible in $\mathcal R^+$ and $h(t^{-1})\in1+t^{-1}\mathcal R^{-}$ invertible in $\mathcal R^-$ 
such that
\[
f= c\cdot t^n \cdot g(t) \cdot h(t^{-1})\,.
\]

Then 
we have:\vspace{3mm}

(4.1)\,{\bf Proposition} {\it
For $f_1,f_2\in \mathscr O_{\widehat 0}^{\ast}$ (or for $f_1,f_2\in \mathcal R^{\ast}$), with decompositions
\[
 f_i= c_i\cdot t^{n_i} \cdot g_i(t) \cdot h_i(t^{-1})
\]
as above ($i=1,2$), we have $(f_1,f_2)_+=c_1^{-n_2}\cdot c_2^{n_1}=(f_1,f_2)^{-1}_-$. }\vspace{2mm}

\begin{Proof} 
We prove the formula in the complex case for $(\cdot, \cdot)_+$, the proofs in the $p$-adic case or for $(\cdot, \cdot)_-$ are analogous.
To lighten notations, set $V=\mathscr O_{\widehat 0}$, $P =\mathscr O_{\widehat 0}^+$,
$N =\mathscr O_{\widehat 0}^-$, the image of $x\in V$ by the projection 
$V=P\oplus N \twoheadrightarrow P$ will be denoted $x_+$. By i) and ii) in (3.4),
it suffices to prove the following statements (cf. \cite{Del}*{2.8.ii}):
\begin{enumerate}
\item[a)] $(g_i,g_j)_+=(g_i,h_j)_+=(h_i,h_j)_+=1$ for $i\neq j\in\{1,2\}$: 
It is enough to prove that the lines $\mathbf P_{g_i}, \mathbf P_{h_i}$ 
are canonically trivial. For $\mathbf P_{g_i}$ this is immediate because $g_i(P)=P$. For $\mathbf P_{h_i}$,
consider the composition $\alpha: P \hookrightarrow V \twoheadrightarrow h_i(P)$, where the second arrow is the 
projection corresponding to the decomposition $V=h_i(P) \oplus h_i(N)$. This map is given by
$\alpha(x)=h_i(h_i^{-1}(x)_+)$. Since $h_i^{-1}$ only involves negative powers of the parameter $t$, the map
\begin{eqnarray*}
 P &\longrightarrow & P \\
x &\longmapsto & h_i^{-1}(x)_+
\end{eqnarray*}
is a bijection, thus 
$\alpha$ is bijective as well and the triviality of $\mathbf P_{h_i}$ follows.

\item[b)] $(c_it^{n_i},g_j)_+=(c_it^{n_i},h_j)_+=1$ for $i\neq j\in\{1,2\}$: 
By bimultiplicativity and the previous case, we can assume that $n_1=n_2=1$ and $c_1=c_2=1$. We will prove that 
$(t, g)_+=1$ if $g$ is holomorphic at zero and $g(0)=1$, 
the other cases are similar.
As just seen, $\mathbf P_{g}\cong \mathbb C$ 
canonically. The projection map $P \longrightarrow tP\,$ is surjective with kernel $\langle 1 \rangle \subset P$ thus 
$\mathbf P_t=\mathbb C\cdot 1^{\ast}$. Take the liftings $\widetilde{t}=(t,1^{\ast})$ and $\widetilde{g}=(g,1)$, 
the definition of the operation on $GL^+(V)$ shows that $\widetilde{t}\cdot \widetilde{g}=(t\cdot g, 1^{\ast})$.
Similarly, $\mathbf P_{g^{-1}}\cong \mathbb C$,
$\mathbf P_{t^{-1}}=\mathbb C \cdot t^{-1} $ and $\widetilde{t^{-1}}\cdot \widetilde{g^{-1}}
=(t^{-1}\cdot g^{-1}, t^{-1})$. Again by definition of the operation on $GL^+(V)$, one easily checks that
\[
\widetilde{t} \ \ \widetilde{g}\ \ \widetilde{t^{-1}}\ \ \widetilde{g}^{-1}=(1,1) \ ,
\]
and so $(t, g)_+=1$.
\item[c)] $(c_1t^{n_1}, c_2t^{n_2})_+=c_1^{-n_2}\cdot c_2^{n_1}$: It is enough to consider the case
$n_1=0, n_2=1, c_2=1$. In this case, a computation as in b)  proves the desired equality. \ \ \ \ $\Box$
\end{enumerate}
\end{Proof}

 For $f \in  \mathscr O_{\widehat 0}$, 
let $\widetilde{f}$ be a representative of $f$ defined on a small circle $S\subset \mathbb C$ 
around $0$ and denote $v(f)$ the degree of the map $arg(\widetilde{f}): S \longrightarrow \mathbb S^1$. 
In \cite{Del}, a symbol 
$
 (\cdot,\cdot)_{ D}: \mathscr O_{\widehat 0}^{\ast} \times \mathscr O_{\widehat 0}^{\ast} \longrightarrow \mathbb C^{\ast}
$
is defined by: 
\[
 (f, g) _{ D}=\exp\left( -\frac{1}{2\pi i} \int_S \log \widetilde{f} \cdot \widetilde{g}^{-1}d\,\widetilde{g} \right)
\cdot \widetilde{g}(0)^{v(f)}\ .
\]
In case $f,g$ are meromorphic at zero, it is proved in loc. cit. that 
$ (f, g) _{D}
$ equals Tate's tame symbol, that is, $ (f, g) _{D}= (-1)^{v(f)\, v(g)}\cdot 
\left(\frac{g(x)^{v(f)}}{f(x)^{v(g)}}\right)(0)$. It follows from (4.1) that we have:
\vspace{2mm}

\begin{Corollary}
If $f_1,f_2\in \mathscr O_{\widehat 0}^{\ast}$ are meromorphic at zero, then 
$(f,g)_-$ equals the tame symbol.
\end{Corollary}

\begin{Remark}
It is easy to see that if the functions $f,g$ are not meromorphic, then 
the symbol $(f,g)_-$ does not need to coincide with Deligne's symbol.  
\end{Remark}

We prove next two reciprocity laws for the symbols defined above. We will need the following lemma:
\vs

{\bf Lemma:}
{\it Let $V=V_+\oplus V_-$ be a polarized vector space. Assume we have a decomposition
$V=V^0\oplus V^1$, put $V_i^j=V_i\cap V^j$ for $i\in\{+,-\}$ and $j\in\{0,1\}$. 
If $V^0=V_+^0\oplus V_-^0$ and $V^1=V_+^1\oplus V_-^1$, then $V^0$ and $V^1$ are polarized vector spaces, 
$ GL^+_{res}(V^0) \times GL^+_{res}(V^1)$ can be identified with a subgroup of $GL^+_{res}(V)$, 
and given commuting elements
\begin{eqnarray*}
  f&=&(f_0,f_1)\in GL^+_{res}(V^0) \times GL^+_{res}(V^1) \\
 g&=&(g_0,g_1)\in GL^+_{res}(V^0) \times GL^+_{res}(V^1)
\end{eqnarray*}
we have
\[
 (f,g)_+ = (f_0,g_0)_+ \cdot (f_1, g_1)_+.
\]
A corresponding statement holds for $(\cdot, \cdot)_-$.}
\begin{Proof} It follows from elementary properties of determinants that we have natural isomorphisms $\mu_f: \mathbf P_{f_0} \otimes \mathbf P_{f_1} \longrightarrow \mathbf P_f $
and $\mu_g: \mathbf P_{g_0} \otimes \mathbf P_{g_1} \longrightarrow \mathbf P_g$. Choose liftings 
$\widetilde{f_0}=(f_0,\alpha_{f_0})$, $\widetilde{f_1}=(f_1,\alpha_{f_1})$, 
$\widetilde{g_0}=(g_0,\alpha_{g_0})$, $\widetilde{g_1}=(g_1,\alpha_{g_1})$.  \vspace{1mm}

Then $\tilde{f}=((f_0,f_1), \mu_f(\alpha_{f_0}\otimes \alpha_{f_1}))$ and $\tilde{g}=((g_0,g_1), \mu_g(\alpha_{g_0}\otimes \alpha_{g_1}))$
are liftings to $GL^+_{res}(V)$ of $f$ and $g$. Computing the commutator
$\tilde{f}\,\tilde{g}\,\tilde{f}^{-1}\tilde{g}^{-1}$ with these liftings one 
obtains the desired formula, the details are left to the reader.
$\Box$
\end{Proof}

\begin{Theorem}
Let $X/\mathbb C$ be a complete non-singular curve, $S\subset X$ a finite subset. For $s\in S$, let 
$\mathscr O_{X,\widehat{s}}$
denote the ring of function germs which are holomorphic in a punctured neighborhood of $s$ in $X$.
If $f,g$ are invertible elements of the ring of analytic functions on $X-S$, then
\[
 \prod_{s\in S}(f_s,g_s)_{s}=1.
\]
where $(\cdot,\cdot)_s$ denotes either the minus or the plus symbol in $\mathscr O_{X,\widehat{s}}$.
\end{Theorem}

\begin{Proof} It is enough to consider the minus case. 
We claim first that $
V=\bigoplus_{s\in S} \, \mathscr O_{X,\widehat{s}}
$ is an analytic Tate space and $V_-:=Im[H^0(X-S, \mathscr O_X)\longrightarrow  \bigoplus_{s\in S} \mathscr O_{X,\hat s}]$ 
is a (FN)-lattice in it, the argument is very similar to that in \cite{Tate}*{pg. 156}:
Consider the  sheaves of topological vector spaces on $X$ defined by 
\begin{eqnarray*}
U &\longmapsto& \mathcal F^0(U):=  \mathscr O_X(U-S) \oplus\left(\bigoplus_{s\in U\cap S}\mathscr O^+_{X,\hat{s}}\right) \\
U &\longmapsto& \mathcal F^1(U):= \bigoplus_{s\in U\cap S}\mathscr O_{X,\hat{s}}\, ,  
\end{eqnarray*}
where $U\subset X$ is open and a direct sum is equal to zero if its index set is empty. We have an exact sequence
\[
 0 \longrightarrow \mathscr O_X \longrightarrow \mathcal F^0 \stackrel{\delta}\longrightarrow \mathcal F^1 
\longrightarrow 0\, ,
\]
where $\delta(\,f\,,\, \oplus_s\,g_s\,)=\oplus_s(f-g_s)_s$, taking global sections we get
\[
0\rightarrow H^0(X,\mathscr O_X)\rightarrow H^0(X-S, \mathscr O_X)\oplus\left(\bigoplus_{s\in S}\mathscr O^+_{X,\hat{s}}\right) 
\stackrel{H^0(\delta)} \rightarrow \bigoplus_{s\in S}\mathscr O_{X,\hat{s}}  
\rightarrow H^1(X, \mathscr O_X).
\]
This may be regarded as a sequence of topological vector spaces (considering \u{C}ech cohomology, see \cite{CaSe}).
The space $H^1(X, \mathscr O_X)$ is Hausdorff [loc.cit.],
and it follows that $H^0(\delta)$ has closed image. Since $H^0(X-S, \mathscr O_X)$ is a (FN)-space,
$\bigoplus_{s\in S}\mathscr O^+_{X,\hat{s}}$ is (DFN), and the open mapping theorem holds for the spaces involved, 
we conclude that $V$ is an analytic Tate space and  $V_-$ is a (FN)-lattice on it, as claimed.
\vs

The theorem can now be proved as in \cite{ADK} or \cite{AR}, considering in the analytic Tate space $V$
the (FN)-lattices $V_-$ and $V_-'=\bigoplus_{s\in S}\, \mathscr O^-_{X,\widehat{s}}$.
They determine the same polarization by corollary (3.6.i), so the value of $(f,g)_-$ is the same if computed 
with $V_-$ or with $V'_-$.
We have $f(V_-)=V_-$ and $g(V_-)=V_-$, so calculating $(f,g)_-$ with $V_-$ we obtain $(f,g)_-=1$. 
To compute the pairing with $V'_-$, use the previous lemma, the claimed equality follows.
$\Box$
\end{Proof}

\vspace{3mm}

Let $\mathbb F$ be a finite field, 
$X/\mathbb F$ a projective non-singular curve, let $Y\subset X$ be an affine curve. Set
\[
A^{\dagger}_Y= \varinjlim_U \Gamma(U, \mathscr O_U)
\]
where $U$ runs over the set of strict neighborhoods of the tube $]Y[$ (we refer 
for these notions to \cite{Cr0}*{section 7}). Let $S$ be the set of closed points of $X-Y$,
we have:

\begin{Theorem}
Under the above assumptions, 
if $f,g$ are invertible elements of $A^{\dagger}_Y$, then
\[
\prod_{s\in S}(f_s,g_s)_s=1.
\]
where $(\cdot,\cdot)_s$ denotes either the minus or the plus symbol calculated in the Robba ring at $s\in S$.
\end{Theorem}
\begin{Proof}
Let $\mathcal R_s$ denote the Robba ring at $s\in S$. By \cite{Cr}*{Theorem 7.5}, $A^{\dagger}_Y$ is a (DFN)-lattice 
in the analytic Tate space $\bigoplus_{s\in S}\mathcal R_s$ (nuclearity is not 
explicitly mentioned in loc.cit., but it can be easily checked). Now the proof is as in the previous theorem, 
considering the (DFN)-lattices $A^{\dagger}_Y$ and $\bigoplus_{s\in S}\mathcal R^{-}_s$.
$\Box$
\end{Proof}

\end{document}